\documentclass[12pt]{article}
\usepackage{amsmath}
\usepackage{amssymb}
\usepackage{amsthm}
\usepackage[totalwidth=17cm, totalheight=25cm]{geometry}
\usepackage{graphicx}
\usepackage{mathabx}
\usepackage{tabularx}
	\newcolumntype{C}[1]{>{\centering\arraybackslash}m{#1}} 
	\newcolumntype{R}[1]{>{\raggedleft\arraybackslash}m{#1}} 

\newtheoremstyle{boldplain}
{9pt}
{9pt}
{\itshape}
{}
{\bfseries}
{.}
{.5em}
{\thmname{#1}\thmnumber{ #2}\thmnote{ (#3)}}%

\newtheoremstyle{bolddefinition}
{9pt}
{9pt}
{}
{}
{\bfseries}
{.}
{.5em}
{\thmname{#1}\thmnumber{ #2}\thmnote{ (#3)}}%

\theoremstyle{boldplain}

\newtheorem{cor}[equation]{Corollary}

\newtheorem{thm}[equation]{Theorem}

\theoremstyle{bolddefinition}
\newtheorem{dfn}[equation]{Definition}

\newfont{\bigbf}{cmbx10 scaled\magstep1}
\normalbaselines
\numberwithin{equation}{section}

\def\no{\noindent}

\def\R{{\mathbb R}}

\def\N{{\mathbb N}}

\def\Z{{\mathbb Z}}

\def\al{\alpha}
\def\ga{\gamma}
\def\Ga{\Gamma}

\def\De{\Delta}
\def\eps{\epsilon}

\def\La{\Lambda}
\def\si{\sigma}
\def\Si{\Sigma}

\def\Om{\Omega}

\def\3{\ss}

\def\acts{\curvearrowright}
\def\amod{a_{mod}}

\def\barXt{\operatorname{\bar X}^{\taumod}}

\def\D{\partial}
\def\DF{\partial_{F\ddot u}}
\def\Ds{\partial_{\si_{mod}}}
\def\Dt{\partial_{\tau_{mod}}}

\def\diamos{\diamondsuit_{\si_{mod}}}
\def\diamot{\diamondsuit_{\tau_{mod}}}

\def\embed{\hookrightarrow}

\def\Flagt{\operatorname{Flag_{\tau_{mod}}}}

\def\geo{\partial_{\infty}}

\def\Lat{\La_{\tau_{mod}}}

\def\lra{\longrightarrow}

\def\ol{\overline}
\def\OmF{\Om_{F\ddot u}}

\def\2pithird{\frac{2\pi}{3}}

\def\rank{\mathop{\hbox{rank}}}
\def\Rep{\operatorname{Rep}}

\def\Ra{\Rightarrow}

\def\simod{\si_{mod}}

\def\st{\operatorname{st}}

\def\taumod{\tau_{mod}}

\def\Th{\mathop{\hbox{Th}}\nolimits}
\def\ThF{\mathop{\hbox{Th}}_{F\ddot u}\nolimits}
\def\ThFi{\mathop{\hbox{Th}}^{Fins}\nolimits}

\def\8{\infty}
\def\<{\langle}
\def\>{\rangle}

\def\BI{\begin{itemize}}
\def\EI{\end{itemize}}

\long\def\comment#1\endcomment{}

\title{Some recent results on Anosov representations}
\author{Michael Kapovich, Bernhard Leeb, Joan Porti}
\date{November 28, 2015}

\begin{document}

\maketitle

\begin{abstract}
In this note 
we give an overview of
some of our recent work on Anosov representations of discrete groups into higher rank semisimple 
Lie groups. 
\end{abstract}

The subject of Anosov representations 
grows out of classical Kleinian groups theory and higher Teichm\"uller theory
originated by Hitchin in \cite{Hitchin}.
Anosov representations of surface groups were introduced by Labourie 
in the rich and beautiful
paper \cite{Labourie}
where he studies the geometry of the representations $\pi_1(\Si)\to PSL(n,\R)$ 
of fundamental groups of compact hyperbolic surfaces $\Si$
which are contained in the Hitchin components of the representation variety $\Rep(\pi_1(\Si),PSL(n,\R))$.
The notion of Anosov representation was subsequently extended 
to representations 
from word hyperbolic groups into semisimple Lie groups 
in \cite{GW}.

The goal of this note is to give an overview and a unified discussion 
of some of the main results
of our papers
\cite{coco13,morse,mlem,bordif2,coco15}. 
We first present a 
``flow-free" definition of Anosov subgroups
and discuss various equivalent dynamical and geometric characterizations of them,
generalizing (to higher rank) characterizations of convex cocompactness in the theory of Kleinian groups.
We put particular emphasis on describing the (coarse extrinsic) geometry of Anosov subgroups,
notably the Morse property.
The coarse geometric viewpoint also leads to 
our local-to-global principle for the Anosov property 
and its application to the construction of Anosov Schottky subgroups. 
Afterwards
we explain our results on the topological dynamics (domains of proper discontinuity and cocompactness) 
of discrete group actions on flag manifolds and Finsler compactifications of symmetric spaces of noncompact type,
and their application to the construction of bordifications and compactifications 
of certain locally symmetric spaces of infinite volume. 

\setcounter{section}{1}
\setcounter{equation}{0}

{\bf 1.\ Definition.}
Let us begin by explaining what Anosov representations are. 
We give here our version
of the original definition in \cite[\S 2-3]{Labourie} and \cite[\S 2]{GW},
where the geodesic flow of the hyperbolic group is replaced by a coarse-geometric object,
namely by the family of (coarse) geodesic lines in the group with respect to a word metric. 
Our definition is equivalent
and has the advantage that it makes the notion of Anosov subgroup
technically more accessible and easier to work with, 
see \cite[\S 6.5.1]{morse} for a detailed discussion and comparison of both definitions.
A (quasi-)geodesic flow on word hyperbolic groups 
had been constructed originally by Gromov \cite{Gromov_hypgps} and later improved by 
Champetier \cite{Champetier} and Mineyev \cite{Mineyev}.
While it has a nice realization e.g.\ for fundamental groups of closed manifolds of negative sectional curvature,
it is in general a technically quite involved object. 

Throughout this note,
let $G$ be a noncompact connected semisimple real algebraic group.
In a nutshell,
a representation $\rho:\Ga\to G$ of a word hyperbolic group $\Ga$ into $G$ is Anosov
if the induced action on a flag manifold associated to $G$ 
satisfies an asymptotic expansion property relative to a continuous boundary map into the flag manifold.

We first describe more precisely what kind of boundary maps are considered
and then state the expansion condition. 

Recall that 
the conjugacy classes of parabolic subgroups $P<G$ 
one-to-one correspond to the faces $\taumod$ of the spherical Weyl chamber $\simod$ attached to $G$.
The conjugacy classes 
yield natural compact homogeneous $G$-spaces (actually, smooth projective varieties) 
$$ \Flagt \cong G/P $$
called {\em (generalized partial) flag manifolds},
and may be viewed as partial {\em boundaries} attached to $G$ 
at infinity.
For $G=SL(n,\R)$, these are precisely the partial flag manifolds.

Two parabolic subgroups are said to be {\em opposite} or {\em antipodal} 
if they are the stabilizers of opposite simplices in the spherical Tits building associated to $G$,
equivalently, 
if they can be swapped by a Cartan involution of $G$.

We will assume in the sequel that $\taumod$ is fixed by the opposition involution 
of $\simod$.
Then 
opposite parabolic subgroups are conjugate to each other. 
We call two {\em flags}, that is, points in $\Flagt$ {\em opposite}
if their parabolic stabilizers in $G$ are.
\begin{dfn}[Boundary embedded (cf.\ {\cite[Def.\ 6.18]{morse}})]
\label{def:anos}
We say that a representation $\rho:\Ga\to G$ 
of a word hyperbolic group $\Ga$ 
is {\em $\taumod$-boundary embedded} 
if there exists a $\Ga$-equivariant topological embedding  
\begin{equation*}
\beta: \geo \Ga \embed \Flagt
\end{equation*}
which is {\em antipodal} in the sense that it
maps different boundary points to opposite flags.
\end{dfn}
Note that boundary embedded representations are 
discrete and have finite kernel,
because $\Ga$ acts on $\beta(\geo\Ga)\cong\geo\Ga$ as a discrete convergence group.

To make the expansion condition precise,
we equip the flag manifolds with (arbitrary) Riemannian metrics. 
For an element $\ga\in\Ga$,
we measure the expansion of the diffeomorphism $\rho(\ga):\Flagt\to\Flagt$ at a point $\tau\in\Flagt$
in terms of its differential by
\begin{equation*}
\eps(\ga,\tau) := 
\bigl\| \bigl( d(\rho(\ga))_{\tau}\bigr)^{-1} \bigr\|^{-1} .
\end{equation*}
Furthermore, we fix a word metric on $\Ga$ 
and consider the asymptotics of the expansion rates 
for sequences of group elements following coarse geodesic rays in $\Ga$.
(The geodesic rays replace the trajectories of the geodesic flow of $\Ga$ in the original Anosov definition.)
We recall that in a word hyperbolic group geodesic rays converge at infinity in the Gromov compactification. 

\begin{dfn}[Anosov, cf.\ {\cite[Def.\ 6.45]{morse}}]
\label{dfn:our anosov}
We say that the representation $\rho$ is {\em $\taumod$-Anosov}
if it is $\taumod$-boundary embedded 
and if for every ideal point $\zeta\in \geo \Ga$
and every normalized 
(by $q(0)=e\in \Ga$)
geodesic ray $q: \N\to \Ga$ asymptotic to $\zeta$
it holds that 
\begin{equation*}
\eps(q(n)^{-1}, \beta(\zeta))\ge A e^{Cn}
\end{equation*}
for $n\geq 0$ with constants $A, C>0$  independent of $q$. 
\end{dfn}

Our notion of $\taumod$-Anosov is equivalent to the notion of $P$-Anosov in \cite{GW}
where $P<G$ is a parabolic subgroup in the conjugacy class corresponding to $\taumod$.
Note also that the study of $(P_+,P_-)$-Anosov representations quickly reduces 
to the case of $P$-Anosov representations by intersecting parabolic subgroups,
cf.\ \cite[Lemma 3.18]{GW}.
Being boundary embedded, 
Anosov representations have discrete image and finite kernel,
and we will refer to $\rho(\Ga)$ as a {\em $\taumod$-Anosov subgroup} of $G$. 

In both our and the original definition uniform exponential expansion rates are required.
It turns out that the conditions can be relaxed without altering the class of representations.
Uniformity can be dropped, and instead of exponential divergence 
the mere divergence of the expansion rate along a subsequence suffices:
\begin{dfn}
[Non-uniformly Anosov, cf.\ {\cite[Def.\ 6.46]{morse}}] 
We say that the representation $\rho$ is {\em non-uniformly $\taumod$-Anosov} 
if it is $\taumod$-boundary embedded 
and if for every ideal point $\zeta\in \geo \Ga$
and every normalized\footnote{Here, the normalization can be dropped because no {\em uniform} growth is required.}
geodesic ray $q: \N\to \Ga$ asymptotic to $\zeta$
it holds that 
\begin{equation*}
\limsup_{n\to+\infty}\, \eps(q(n)^{-1}, \beta(\zeta))=+\infty. 
\end{equation*}
\end{dfn}
In other words,
we require that for every ideal point $\zeta\in\geo\Ga$
the expansion rate $\eps(\ga_n^{-1},\beta(\zeta))$
non-uniformly diverges along some sequence $(\ga_n)$ in $\Ga$
which converges to $\zeta$ {\em conically}.

It is shown in \cite[Thm.\ 6.57]{morse} that 
non-uniformly $\taumod$-Anosov implies $\taumod$-Anosov. 
Other relaxations of the original Anosov condition also appear in \cite[sec.\ 6.1]{Labourie} and \cite[Prop.\ 3.16]{GW}.
They assume uniform, not necessarily exponential, divergence of the expansion factors  of (lifted) geodesic flows. 

\setcounter{section}{2}
\setcounter{equation}{0}

\medskip
{\bf 2.\ Equivalent characterizations.}
We will now discuss
various dynamical and geometric characterizations of Anosov subgroups.

In the case of discrete subgroups of {\em rank one} Lie groups, 
e.g.\ for Kleinian groups,
the Anosov condition is equivalent to {\em convex cocompactness}. 
We refer the reader to 
the paper 
\cite{Bowditch_gf} for a list of equivalent characterizations of convex cocompact isometry groups of negatively curved spaces. 
Whereas in higher rank convex cocompactness is too restrictive, 
compare \cite{convcoco,Quint},
we will see that suitable versions of other equivalent conditions considered in the theory of Kleinian groups 
remain meaningful 
and provide alternative characterizations for Anosov subgroups.
Having a supply of non-obviously equivalent 
characterizations 
makes it possible to switch back and forth between different viewpoints,
geometric and dynamical,
in a fruitful way.

{\bf 2.1.\ Dynamical characterizations.}
Let us first remain 
in the framework of 
dynamics on flag manifolds 
and explain a certain {\em convergence type dynamics}
enjoyed by Anosov subgroups, 
cf.\ \cite[\S 6]{coco15},
and how Anosov subgroups are distinguished among subgroups with this kind of dynamics.

For a flag $\tau\in\Flagt$ let us denote by $C(\tau)\subset\Flagt$
the {\em open Schubert stratum} consisting of the flags opposite to $\tau$.
It is an open and dense orbit of the parabolic stabilizer $P_{\tau}$ of $\tau$,
and its complement is a projective subvariety.

We say 
that a sequence $(g_n)$ in $G$ is {\em $\taumod$-contracting}
if there exist flags $\tau_{\pm}\in\Flagt$ such that 
\begin{equation} 
\label{eq:convoncell}
g_n|_{C(\tau_-)}\to\tau_+
\end{equation}
uniformly on compacta as $n\to+\infty$. 

\begin{dfn}
[Convergence]
A discrete subgroup $\Ga<G$ is called a {\em $\taumod$-convergence subgroup}
if every sequence $\ga_n\to\infty$ in $\Ga$ contains a $\taumod$-contracting subsequence.
\end{dfn}

Note that in rank one
this coincides with the usual notion of {\em convergence subgroup} 
\cite{Bowditch_char,Bowditch_config}
and is satisfied by all discrete subgroups.
(In rank one, 
the open Schubert strata are the complements of points in the visual boundary.)

For an arbitrary discrete subgroup $\Ga<G$ we define the {\em $\taumod$-limit set}
$$ \Lat \subset\Flagt$$
as the set of all flags $\tau_+$ as in (\ref{eq:convoncell}) for all $\taumod$-contracting sequences in $\Ga$.
It is compact and $\Ga$-invariant.
The structure of the dynamics $\Ga\acts\Flagt$ is particularly closely related to the limit set 
in the case of $\taumod$-convergence subgroups. 
We refer the reader to \cite{Benoist} for discussion 
of related notions of limit sets of discrete subgroups of semisimple Lie groups and their properties.

Anosov subgroups are $\taumod$-convergence subgroups and, among these,
can be distinguished in different ways.
One possibility is a stronger form of boundary embeddedness
where one requires convergence dynamics 
and that the boundary map identifies the Gromov boundary of the discrete subgroup 
with its limit set in the flag manifold. 
Here, we call a subset of $\Flagt$ {\em antipodal} if it consists of pairwise opposite flags. 

\begin{dfn}[Asymptotically embedded, cf.\ {\cite[Def.\ 1.5]{morse}}]
We say that a $\tau_{mod}$-convergence subgroup $\Ga<G$ 
is {\em $\taumod$-asymptotically embedded}, 
if 

(i) $\Lat$ is antipodal, 

(ii) $\Ga$ is intrinsically word hyperbolic 
and there exists a $\Ga$-equivariant homeomorphism 
\begin{equation*}
\label{eq:mapalphatauintro}
\al:\geo\Ga \stackrel{\cong}{\lra} \Lat .
\end{equation*}
\end{dfn}

Another possibility to distinguish Anosov among convergence subgroups 
is by an {\em expansivity property} 
which was originally introduced by Sullivan 
in his study of Kleinian groups.
We formulate it in a general setting:

\begin{dfn}[Expanding action (cf.\ {\cite[\S 9]{Sullivan}})]
\label{dfn:expact}
A continuous action $\Ga\acts Z$ of a discrete group $\Ga$ on a compact metric space $(Z,d)$
is said to be 
{\em expanding} 
at the {\em point} $z\in Z$ 
if there exists an element $\ga\in\Ga$ 
which is {\em uniformly expanding} on a neighborhood $U$ of $z$, 
i.e.\ for some constant $c>1$ and all points $z_1,z_2\in U$ we have  
\begin{equation*}
d(\ga z_1,\ga z_2)\geq c\cdot d(z_1,z_2) .
\end{equation*}
The action is said to be {\em expanding} 
at a compact $\Ga$-invariant {\em subset} $E\subset Z$ 
if it is expanding at all points $z\in E$. 
\end{dfn}

Returning to our framework,
we define the following class of subgroups.
\begin{dfn}
[Expanding, cf.\ {\cite[Def 7.12]{coco15}}]
We say that a $\taumod$-convergence subgroup $\Ga<G$ 
is {\em $\taumod$-CEA}
(convergence, expanding, antipodal)
if $\Lat$ is antipodal 
and if the action $\Ga\acts\Flagt$ is expanding at $\Lat$.
\end{dfn}
The following result is obtained by combining \cite[Thms.\ 1.7 and 5.23]{morse}:

\begin{thm}[Dynamical characterizations of Anosov I: actions on flag manifolds]
For a discrete subgroup $\Ga<G$,
the following properties are equivalent:

(i) $\taumod$-Anosov, 

(ii) $\taumod$-CEA,

(iii) $\taumod$-asymptotically embedded. 
\end{thm}

Note that 
CEA 
does not a priori assume word hyperbolicity of 
$\Ga$.
This is a consequence.

\setcounter{section}{3}
\setcounter{equation}{0}

\medskip
{\bf 3.\ Geometric characterizations.}
We will now discuss other characterizations of Anosov subgroups 
which involve 
the geometry of the associated symmetric space of noncompact type $X=G/K$.
Here, $K$ is a maximal compact subgroup of $G$.

A representation $\rho:\Ga\to G$ of an arbitrary discrete group $\Ga$ corresponds to an isometric action
$$\Ga\acts X.$$ 
The action is properly discontinuous iff the representation is discrete with finite kernel.
We will henceforth identify $\Ga$ with its image in $G$ and assume that $\Ga<G$ is a discrete subgroup.

{\bf 3.1.\ Regularity.}
Let us first explain how the dynamical $\taumod$-convergence property can be reformulated
as a regularity condition on the asymptotic geometry of the orbits $\Ga x\subset X$. 

Consider the set 
\begin{equation*}
\{ d_{\De}(x,\ga x):\ga\in\Ga\} \subset\De
\end{equation*}
of $\De$-distances between orbit points.
Here, 
$d_{\De}$ denotes the 
vector-valued {\em $\De$-distance} on $X$ 
with values in the euclidean Weyl chamber $\De$ of $X$, 
compare \cite{KLM} and \cite{Parreau}.
(Recall that $X\times X/G\cong\De$.) 
The $\De$-distance translates into the {\em Cartan projection} in 
Lie theory.

The subgroup $\Ga$ is called {\em regular}, 
if these $\De$-distances 
drift away from the boundary $\D\De$ of the cone $\De$,
i.e.\ if 
\begin{equation}
\label{eq:reg}
d\bigl(d_{\De}(x,\ga x),\D\De\bigr)\to+\infty
\end{equation} 
as $\ga\to\infty$ in $\Ga$,
and {\em uniformly regular},
if the drift is linear in terms of the distance from the origin,
i.e.\ if 
\begin{equation}
\label{eq:ureg}
d\bigl(d_{\De}(x,\ga x),\D\De\bigr)\geq c\cdot\|d_{\De}(x,\ga x) \|-a
\end{equation}
with constants $c,a>0$

We use a relaxation of this regularity condition with respect to $\taumod$
where we measure the drift away from only part of the boundary of $\De$.
Note that the visual boundary at infinity of the euclidean Weyl chamber
is canonically identified with the spherical Weyl chamber, $\geo\De\cong\simod$.
Let $\Dt\simod\subseteq\D\simod$ denote the part of the boundary of $\simod$
which consists of the union of the closed faces not containing $\taumod$.
Accordingly, let $\Dt\De\subseteq\D\De$ be the cone over $\Dt\simod$.
(For instance, $\Ds\simod=\D\simod$ and $\Ds\De=\D\De$.)
\begin{dfn}[Regular, cf.\ {\cite[\S 5]{morse}}]
We say that the subgroup $\Ga<G$ is {\em $\taumod$-regular}, 
respectively, {\em uniformly $\taumod$-regular}
if the corresponding properties (\ref{eq:reg}), respectively, (\ref{eq:ureg}) hold with $\D\De$ replaced by $\Dt\De$.
\end{dfn}
Note that $\simod$-regularity is the same as regularity,
and the reader may just restrict to this special case for simplicity.

\begin{thm}[{\cite[Thm.\ 5.23]{morse}}]
\label{thm:conveqreg}
$\Ga<G$ is $\taumod$-convergence iff it is $\taumod$-regular.
\end{thm}

In particular, $\taumod$-Anosov subgroups are $\taumod$-regular.
They are in fact {\em uniformly} $\taumod$-regular, cf.\ \cite[Thm.\ 6.33]{morse}.

{\bf 3.2. Coarse extrinsic geometry.}
The Anosov property has strong implications for the coarse extrinsic geometry of a subgroup $\Ga<G$
which we discuss now. 

These implications result from control on the images of coarse geodesic segments (rays, lines) in $\Ga$
under the orbit maps $\Ga\to\Ga x\subset X$.
Already 
$\taumod$-boundary embeddedness 
implies by a simple compactness argument\footnote{cf.\ \cite[\S 6.4.1]{morse}}
that the images of coarse lines $q:\Z\to\Ga$ 
under the orbit maps 
are uniformly close\footnote{in the sense of being contained in tubular neighborhoods with uniform radii}
to $\taumod$-parallel sets 
$P(q)\subset X$ 
(if $\taumod=\simod$, to maximal flats)
which are picked out by the boundary map $\beta$.
The expansion part of the 
$\taumod$-Anosov condition 
then allows\footnote{cf.\ \cite[\S 6.5.2, Lemma 6.54]{morse}} 
to further restrict the position of the image paths 
$qx$ along these parallel sets.
Namely, the images $qx|_{\N_0}$ of coarse rays are uniformly close to $\taumod$-Weyl cones
(if $\taumod=\simod$, to euclidean Weyl chambers)
with tips at the initial point $q(0)x$.
Moreover,
they are forced to have a linear drift out to infinity
and a linear drift away from the boundaries of these Weyl cones.
In particular, 
$\taumod$-Anosov subgroups 
satisfy the following property\footnote{cf.\ also \cite[Thm.\ 6.33]{morse}} 
first proven in 
\cite[Prop.\ 3.16 and Thm.\ 5.3]{GW}:

\begin{dfn}[Undistorted]
We say that a finitely generated subgroup $\Ga<G$ is {\em $\taumod$-URU},
if it is 

(i) uniformly $\taumod$-regular, and 

(ii) undistorted,
i.e.\ the inclusion $\Ga\subset G$, equivalently,
the orbit maps $\Ga\to\Ga x\subset X$ are {\em quasi-isometric embeddings} 
with respect to a word metric on $\Ga$.
\end{dfn}

But this notion does not fully capture the control 
on the geometry of the orbits
provided by the above discussion,
on which we elaborate now a little more.

Finite subsegments of the coarse line $q$ in $\Ga$ are the intersection of two subrays in opposite directions. 
We therefore see that 
the finite image paths 
$qx|_{[a_-,a_+]\cap\Z}$ 
are uniformly close to the intersection of two Weyl cones in the parallel set $P(q)$,
opening up towards opposite directions, 
and with tips $x_{\pm}$ uniformly close to the endpoints $q(a_{\pm})x$ of the path.
We call this intersection a {\em $\taumod$-diamond}\footnote{cf.\ \cite[\S 3.5, Def.\ 3.33 and Lemma 3.34]{mlem}}
and denote it by $\diamot(x_-,x_+)$.
For instance, 
in the case $\taumod=\simod$,
if the segment $x_-x_+$ is ($\simod$-)regular,
then the diamond $\diamos(x_-,x_+)$ lies in the unique maximal flat containing $x_-x_+$
and equals the intersection of the two euclidean Weyl chambers with tips at $x_-$ and $x_+$
which contain $x_-x_+$.
There is also a nice description of diamonds from the Finsler geometry viewpoint:
For $G$-invariant ``polyhedral" Finsler metrics on $X$ 
(cf.\ \cite{bordif2})
geodesic segments are no longer unique,
and the diamonds $\diamot(x_-,x_+)$ can be described 
as the unions of all Finsler geodesic segments $x_-x_+$ for suitable Finsler metrics depending on $\taumod$.

We are led to the following definition\footnote{In \cite[Def.\ 7.14]{morse} we work with the smaller ``uniformly regular"
$\Theta$-diamonds instead of diamonds, which makes uniform regularity implicit in the definition.}
which we paraphrase in a non-technical way:
\begin{dfn}[Morse, cf.\ {\cite[\S 7]{morse}}]
(i) 
A {\em $\taumod$-Morse quasigeodesic} in $X$ is a coarsely uniformly $\taumod$-regular quasigeodesic
such that every sufficiently long subpath of it
is uniformly close to the $\taumod$-Weyl diamond 
with tips at the endpoints of the subpath. 

(ii) 
We say that 
a finitely generated subgroup $\Ga<G$ is {\em $\taumod$-Morse}
if the images of uniform coarse quasigeodesics in $\Ga$ under an orbit map $\Ga\to\Ga x\subset X$ 
are uniform $\taumod$-Morse quasigeodsics in $X$. 
\end{dfn}
According to our 
discussion, $\taumod$-Anosov subgroups are $\taumod$-Morse.

The Morse property implies URU, and a priori it seems strictly stronger.
However, it turns out that, conversely, URU implies Morse.
This is a consequence of the following non-equi\-va\-ri\-ant geometric result,
generalizing to higher rank
the classical Morse Lemma for quasigeodesics in (coarsely) negatively curved spaces, 
see \cite{Mostow,Gromov_hypgps}.
The latter asserts that 
uniform quasigeodesic segments in Gromov hyperbolic geodesic metric spaces 
are uniformly Hausdorff close to geodesic segments
with the same endpoints.
In our version of the Morse Lemma for symmetric spaces of arbitrary rank,
we need regularity of the quasigeodesics (which comes for free in rank one)
and the geodesic segments are replaced by the larger diamonds:

\begin{thm}[Morse Lemma 
({\cite[Thm.\ 1.3]{mlem}})]
\label{thm:mlem}
All coarse\-ly uniformly $\taumod$-regular quasigeodesics are $\taumod$-Morse.
\end{thm}
From the Finsler perspective, one can reformulate this result 
to the effect that, 
for suitable Finsler metrics on $X$ depending on $\taumod$, 
such quasigeodesic segments are {\em uniformly Hausdorff close to some Finsler geodesic segment 
with the same endpoints},
compare above our Finsler redefinition of diamonds.
Accordingly, the Morse subgroup property can be viewed 
as a {\em Finsler quasiconvexity property}.\footnote{Recall that 
a subgroup $\Ga$ of a hyperbolic group $\Ga'$ is called {\em quasiconvex} 
if coarse geodesic segments in $\Ga'$ with endpoints in $\Ga$ are uniformly close to $\Ga$.} 

Thus, URU and Morse are equivalent coarse properties of Anosov subgroups.
We show that the latter are actually characterized by these properties.
Essentially,
undistortedness characterizes Anosov subgroups 
among uniformly regular subgroups:

\begin{thm}[Coarse geometric characterizations of Anosov (\cite{mlem},\cite{bordif2})]
\label{thm:cgeomcharanos}
For a finitely generated discrete subgroup $\Ga<G$,
the following properties are equivalent:

(i) $\taumod$-Anosov 

(ii) $\taumod$-URU

(iii) $\taumod$-Morse

(iv) $\taumod$-Finsler quasiconvex

\no
Furthermore,
they imply that $\Ga$ is a coarse equivariant retract.
\end{thm}

Note again that the properties besides Anosov do not a priori assume word hyperbolicity.

The last mentioned retraction property is the following strengthening of undistortedness:

\begin{dfn}
[Retract]
We say that an undistorted finitely generated subgroup $\Ga<G$ is a {\em coarse (equivariant) retract}
if there exist ($\Ga$-equivariant) coarse Lipschitz retractions $G\to\Ga$, equivalently,
$X\to\Ga x$ onto orbits.
\end{dfn}

The previous theorem is a combination of various results:

Our proof of the implication ``Anosov$\Ra$retract"
is based on our proper discontinuity and cocompactness results 
for the actions of Anosov subgroups on suitable domains 
in Finsler compactifications of $X$, 
cf.\ \cite[Thms.\ 1.7]{bordif2} and Theorem~\ref{thm:pdccfins} below.

Regarding the implication ``URU$\Ra$Anosov":
That URU subgroups are, besides being Morse, also intrinsically word hyperbolic
is implied by part (i) of the following 
non-equivariant coarse geometric result:

\begin{thm}[Hyperbolicity and boundary maps ({\cite[Thm.\ 1.4]{mlem}})]
Let $Z$ be a locally compact (quasi)geodesic metric space,
and suppose that $Q: Z\to X$ is a coarsely uniformly $\taumod$-regular
quasiisometric embedding. 
Then 

(i) $Z$ is Gromov hyperbolic, and 

(ii)
$Q$ extends to a map 
$$\bar Q: \bar Z \to \barXt$$
from the (visual) Gromov compactification $\bar Z=Z\sqcup\geo Z$ 
to the $\taumod$-bordification $\barXt=X\sqcup\Flagt$
which is continuous at $\geo Z$ 
and whose restriction $Q|_{\geo Z}$ is antipodal, 
i.e.\ sends distinct ideal boundary points in $\geo Z$ to antipodal flags in $\Flagt$. 
\end{thm}

Note that the bordification $\barXt$ sits in a suitable Finsler compactification, compare \cite{bordif2}.

Part (ii) of the theorem follows from the Morse Lemma~\ref{thm:mlem}.
Applied to orbit maps, 
it then provides the boundary maps 
$\geo\Ga\to\Flagt$
for which our URU subgroups are asymptotically embedded 
\cite[\S 7.4-5]{morse},
and the Morse property translates into the expansion part of the Anosov property,
compare the proof of \cite[Thm.\ 6.57]{morse},
establishing that ``URU$\Ra$Anosov".

{\bf Local-to-global principle.}
In the course of our study of the Morse property \cite[\S 7]{morse},
we also obtain a {\em local-to-global principle} for representations $\rho:\Ga\to G$ of word hyperbolic groups
to be $\taumod$-Anosov,
reminiscent of Gromov's local-to-global principle for the word hyperbolicity of groups, 
namely by verifying a straightness condition for a sufficiently large finite subset of $\Ga$. 
It implies {\em semidecidability} whether the representation is $\taumod$-Anosov.
It also provides proofs of openness and structural stability of Anosov representations, 
properties previously proven in \cite{GW} using a different approach.

{\bf Examples.}
We apply the local-to-global principle in \cite[\S 7.6]{morse} 
to construct Anosov representations of free groups:
Let $a,b\subset X$ be geodesic lines, 
and let $\al,\beta\in G$ be axial isometries with $a,b$ as axes.
For $m,n\in\N$ we consider the representation 
of the free group in, say, two generators 
$$\rho_{m,n}:F_2=\langle A,B\rangle\to G$$
sending $A\mapsto\al^m$ and $B\mapsto\beta^n$. 

\begin{thm}[Anosov Schottky subgroups, cf.\ {\cite[Thm.\ 7.40]{morse}}]
If $(a,b)$ is a generic pair of $\taumod$-regular geodesic lines,
then for sufficiently large $m,n$ 
(depending also on $\al,\beta$)
the representation $\rho_{m,n}$ is faithful and $\taumod$-Anosov.
\end{thm}
Genericity is defined in terms of the relative position of the quadruple of ideal endpoints of $a,b$ 
in the visual boundary. 
This result has later been reproven by a different method, 
going back to the work of Benoist \cite{Benoist}, in the special case $G=SL(d,\R)$ in 
\cite[Thm.\ A.2]{CLS}. 

\medskip
{\bf 3.3. Asymptotic geometry.}
The next condition on discrete subgroups 
generalizes the {\em conical limit points} condition on the asymptotics of orbits 
from the theory of Kleinian groups.

We say that a limit flag $\tau\in\Lat$ is {\em conical} 
if there exists a sequence $\ga_n\to\infty$ in $\Ga$
such that the corresponding orbit sequence $(\ga_nx)$ in $X$ 
is contained in a tubular neighborhood of the Weyl cone $V(x,\st(\tau))\subset X$, 
cf.\ \cite[Def.\ 5.2]{Albuquerque} and \cite[Def.\ 6.1]{morse}.
In the case $\taumod=\simod$,
for a Weyl chamber $\si\subset\geo X$ in the visual boundary,
the Weyl cone $V(x,\st(\si))$ is simply the euclidean Weyl chamber $V(x,\si)$ 
with tip $x$ and asymptotic to $\si$.\footnote{In general,
$V(x,\st(\tau))$ is obtained by coning off at $x$ the star $\st(\tau)\subset\geo X$
which in turn is defined as the union of all spherical Weyl chambers having the flag $\tau$,
thought of as a simplex in the visual boundary, as a face.}

\begin{dfn}[Conical, cf.\ {\cite[Def.\ 1.3]{morse}}]
We say that a subgroup $\Ga<G$ is {\em $\taumod$-RCA} 
if it is $\taumod$-regular
and if all flags in $\Lat$ are conical and pairwise antipodal.
\end{dfn}

Anosov subgroups are RCA,
as follows from them being Morse and asymptotically embedded,
compare our above discussion of their extrinsic geometry.
The converse holds as well:

\begin{thm}[Asymptotic geometric characterization of Anosov ({\cite[Thm.\ 1.7]{morse}})]
\label{thm:cgeomcharanos2}
A subgroup $\Ga<G$ is $\taumod$-Anosov iff it is $\taumod$-RCA.
\end{thm}

Note again that the RCA property does not a priori assume word hyperbolicity 
or even finite generation of $\Ga$. 

The implication ``RCA$\Ra$Anosov" 
is obtained by observing that, due to antipodality, 
the restricted action $\Ga\acts\Lat$ is a convergence action
and by translating our extrinsic conicality condition into the intrinsic one 
in terms of the action $\Ga\acts\Lat^3$ on triples,
cf.\ \cite[\S 6.1.4]{morse}.
Using Bowditch's characterization \cite{Bowditch_char} 
of hyperbolic groups in terms of their dynamics at infinity,
we then conclude that $\Ga$ is word hyperbolic and asymptotically embedded,
compare the proof of \cite[Thm.\ 6.16]{morse}.

\setcounter{section}{4}
\setcounter{equation}{0}

\medskip
{\bf 4.\ Topological dynamics.}
We explain now our main results regarding the topological dynamics of discrete group actions 
on flag manifolds \cite{coco13,coco15} and Finsler compactifications \cite{bordif2}.

{\bf 4.1.\ Convergence dynamics in rank one and implications.}
Recall that in {\em rank one} (e.g.\ for Kleinian groups)
the action 
$$\Ga\acts\ol X=X\sqcup\geo X$$
on the visual compactification has {\em convergence dynamics}.
This leads to the $\Ga$-invariant {\em dynamical decomposition}
$$ \ol X=\underbrace{X\sqcup\Om_{\infty}}_{\Om}\sqcup\La $$
into the {\em domain of discontinuity} $\Om$, 
obtained from $X$ by attaching the domain of discontinuity {\em at infinity} $\Om_{\infty}$, 
and the {\em limit set} $\La$.
It also yields that the action on $\Om$ is not just discontinuous, but {\em properly} discontinuous.
Furthermore, 
$\Ga$ is convex cocompact iff the action $\Ga\acts\Om$ is cocompact,
see \cite{Bowditch_gf}. 
In this case, in particular the action $\Ga\acts\Om_{\infty}$ at infinity is cocompact.

{\bf 4.2.\ Visual versus Finsler compactifications.}
In rank $\geq2$, 
the situation is more complicated 
as convergence dynamics gets lost (at least in its full strength).
The visual compactification becomes less suitable for finding good domains of discontinuity 
and we work with a different compactification
suggested by a Finsler geometry point of view.
(Note that 
Finsler geometry appeared implicitly 
already in our earlier discussion,
like in the notion of regularity and in the reformulation of the Morse Lemma, 
cf.\ Theorem~\ref{thm:mlem} above.)

Symmetric spaces of noncompact type are CAT(0) spaces.
The visual compactification of proper CAT(0) spaces 
can be seen as a special case 
of a very general procedure for producing geometric compactifications,
the {\em horoboundary} construction, 
which applies to any proper geodesic metric space $Y$ 
(cf.\ \cite{Gromov_hypmfs}, \cite[\S II.1]{Ballmann}).
One considers the natural topological embedding 
\begin{equation*}
Y\lra \ol{\mathcal C}(Y), \quad y\mapsto [d_y] .
\end{equation*}
into the space $\ol{\mathcal C}(Y)={\mathcal C}(Y)/\R$
of continuous real valued functions modulo additive constants,
assigning to a point $y\in Y$ the equivalence class of the distance function $d_y=d(y,\cdot)$.
By taking the closure of the image of $Y$, one obtains the {\em horoclosure}
$$ \ol Y =Y\sqcup\geo Y .$$
Returning to the symmetric space $X$,
if one replaces the {\em Riemannian} metric 
by suitable $G$-invariant polyhedral {\em Finsler} metrics,
other compactifications arise 
with different asymptotic geometry and dynamics at infinity,
and which turn out to be more suitable for our purposes.

A $G$-invariant Finsler metric on $X$ corresponds to a 
Weyl group invariant norm on a Cartan subgroup $A$ of $G$.
It is determined by its unit ball which we choose to be a polyhedron
with vertices on the maximally singular rays in $A$ emanating from the origin;
more precisely, we require that there is one (non-degenerate) vertex on every singular ray
and that there are no other vertices. 
The resulting {\em regular Finsler compactification} $$\ol X^{Fins}=X\sqcup\geo^{Fins}X$$
does not depend on the particular choice of the Finsler metric.
It is a geometric realization of the {\em maximal Satake compactification} from the theory of algebraic groups
and has especially nice topological, geometric and dynamical properties 
(see \cite{bordif2} and \cite{Parreau} for a careful discussion).
The most important ones are:
\begin{thm}[{\cite[Thm.\ 1.1]{bordif2}}]
(i) The natural action $G\acts\ol X^{Fins}$
has finitely many orbits $S_{\taumod}$ 
which correspond to the faces $\taumod$ of $\simod$, including the empty face; 
the orbit closure inclusion ``$\subseteq$" corresponds to containment of faces ``$\supseteq$".

(ii) The $G$-orbits are the strata of a manifold-with-corners\footnote{The local model
for a $d$-manifold with corners is the $d$-orthant $[0,+\infty)^d$.} structure.

(iii) 
$\ol X^{Fins}$ is homeomorphic to a compact ball. 

(iv)
$\ol X^{Fins}$ is $G$-equivariantly homeomorphic to the maximal Satake compactification $\ol X^S_{max}$.
\end{thm}
Let us add a few more details to the picture:
There is exactly one closed $G$-orbit in $\ol X^{Fins}$,
namely the full flag manifold or {\em F\"urstenberg boundary}
$$S_{\simod}\cong \DF X\cong G/B,$$
and, on the opposite extreme,
the orbit $S_{\emptyset}=X$ is open and dense. 
Here $B< G$ is a Borel subgroup.
The orbits at infinity, 
i.e.\ the strata $S_{\taumod}$ for $\taumod\neq\emptyset$,
are {\em blow-ups} of the corresponding flag manifolds $\Flagt$.
More precisely,
there are $G$-equivariant fibrations 
$$S_{\taumod}\lra\Flagt$$
with contractible fibers.
The fiber $S_{\tau}\subset S_{\taumod}$ over $\tau\in\Flagt$ 
can be interpreted geometrically 
as the space of {\em strong asymptote classes of Weyl sectors} $V(x,\tau)$ asymptotic to $\tau$,
cf.\ \cite[\S 3]{bordif2}. 
In particular, it is a symmetric space of rank 
$\dim\simod-\dim\taumod<\rank X=1+\dim\simod$.
We refer to the $S_{\tau}$ as {\em small strata} at infinity. 
Their boundaries $\D S_{\tau}=\ol S_{\tau}-S_{\tau}$ 
are unions of small strata,
namely of the $S_{\nu}$ 
for the flags $\nu$ strictly ``refining" $\tau$
in the sense that $\nu\supsetneq\tau$ for the corresponding simplices in the visual boundary $\geo X$.

We say that a subset of $\geo^{Fins}X$ is {\em saturated} if it is a union of small strata. 

It is worth noting that 
the stabilizers of the points in the Finsler compactification 
are {\em pairwise different} closed subgroups of $G$.
The stabilizers of the points at infinity in $S_{\tau}$ are contained in the parabolic subgroup $P_{\tau}$.

{\bf 4.3.\ Domains of proper discontinuity and cocompactness.}
Now we can state our dynamical results in higher rank. 

There is a well-defined open {\em domain of discontinuity} or {\em wandering set} $\Om_{disc}\subset\ol X^{Fins}$
for the action $\Ga\acts\ol X^{Fins}$, 
including $X$ itself.
(It consists of the points which have neighborhoods $U$ such that $U\cap\ga U\neq\emptyset$
for at most finitely many $\ga\in\Ga$.)
However,
in higher rank, the $\Ga$-action on the domain of discontinuity is no longer proper
and one has to look for {\em domains of proper discontinuity} $\Om\subset\Om_{disc}$.
Moreover, 
there are in general {\em no unique maximal} such domains.

\begin{thm}[Proper discontinuity and cocompactness on domains in Finsler compactifications,
cf.\ {\cite[Thm.\ 1.8]{bordif2}}]
\label{thm:pdccfins}
Let $\Ga<G$ be uniformly $\taumod$-regular. 
Then there exist natural $\Ga$-invariant saturated open subsets
$\Om_{\infty}\subset\geo^{Fins} X$
such that the following holds:

(i) The action 
$$\Ga\acts\Om=X\sqcup\Om_{\infty}\subset\ol X^{Fins}$$
is properly discontinuous.

(ii) If $\Ga$ is $\taumod$-Anosov, then the action is also cocompact. 
\end{thm}

We show furthermore a converse to the cocompactness part (ii),
implying that Anosov subgroups are characterized among uniformly regular subgroups 
by the cocompactness of this action. 
More generally, we consider the following property:

\begin{dfn}
[$S$-cocompact, cf.\ {\cite[Def.\ 11.2]{bordif2}}]
We say that a discrete subgroup $\Ga<G$ is {\em $S$-cocompact} 
if there exists a $\Ga$-invariant saturated open subset $\Om_{\infty}\subset\geo^{Fins}X$ 
such that the  action 
$\Ga\acts X\sqcup\Om_{\infty}$ 
is properly discontinuous and cocompact. 
\end{dfn}

\begin{thm}[Cocompactness implies Anosov, cf.\ {\cite[Thm.\ 1.9]{bordif2}}]
Uniformly $\taumod$-re\-gu\-lar $S$-cocompact subgroups $\Ga<G$ are are $\taumod$-Anosov.
\end{thm}

We conclude:

\begin{cor}
[Dynamical characterizations of Anosov II: actions on Finsler compactifications,
cf.\ {\cite[Cor.\ 1.10]{bordif2}}]
For a uniformly $\taumod$-regular subgroup $\Ga<G$,
the following properties are equivalent:

(i) $\taumod$-Anosov, 

(ii) $S$-cocompact. 
\end{cor}

{\bf 4.4.\ Locally symmetric spaces.}
Theorem~\ref{thm:pdccfins} yields
bordifications and, in the Anosov case, compactifications 
of the locally symmetric spaces $X/\Ga$ of infinite volume. 
\begin{cor}[Bordifications and compactifications]
Let $\Ga<G$ be uniformly $\taumod$-regular. 
Then there exist natural (real analytic) bordifications 
$$ (X\sqcup\Om_{\infty})/\Ga$$
of the locally symmetric space $X/\Ga$
as orbifolds with corners.

If $\Ga$ is $\taumod$-Anosov, then these bordifications are compactifications.
\end{cor}

The real analyticity comes from the fact that the maximal Satake compactification 
is known to carry a $G$-invariant real-analytic structure,
see \cite{BorelJi}.

{\bf 4.5.\ Thickenings and GIT.}
The domains appearing in Theorem~\ref{thm:pdccfins} arise from a natural construction,
which we can only sketch here. 
We think of their complements
$$ \ThFi(\Lat) = \geo^{Fins}X - \Om_{\infty} $$
as {\em thickenings} of the limit set $\Lat$.
(Compare above the dynamical decomposition of the action $\Ga\acts\ol X$ in the rank one case).
These thickenings admit natural fibrations
$\Th^{Fins}(\Lat)\to\Lat$
with compact saturated fibers.
They 
are derived from a combinatorial datum
which gives a certain degree of flexibility to the construction,
namely from a {\em thickening}
$$ \Th\subset W$$
of the neutral element $e$ in the Weyl group $W$.
By this we mean a union of sublevels for the (strong) Bruhat order, 
cf.\  \cite[\S 8.4]{bordif2} and \cite[\S 3.4]{coco15}.

Suppose that $\Th$ is invariant under the stabilizer $W_{\taumod}<W$ of $\taumod$,
where we think of $\simod$ as embedded in the spherical Coxeter complex $\amod$.
Then it induces a thickening $$\ThF(\Lat)\subset\DF X$$ of 
$\Lat\subset\Flagt$
by taking for every flag $\tau\in A$, thought of as a simplex $\tau\subset\geo X$, 
all chambers $\si\subset\geo X$ 
which have the same position relative to $\tau$ 
as the chambers $w\simod\subset\amod$ for the elements $w\in\Th$ have relative to $\taumod$.
The thickening $\ThFi(\Lat)$ in turn is obtained from $\ThF(\Lat)$ by some filling in procedure. 

The conclusion of Theorem~\ref{thm:pdccfins} holds 
if the thickening $\Th$ is {\em balanced} in the sense that there is the partition 
$$ W =\Th\sqcup w_0\Th$$
where $w_0\in W$ denotes the longest element.

In \cite{coco13,coco15},
we prove more general dynamical results for discrete group actions on {\em flag manifolds}.
For instance, 
for the action 
$$\Ga\acts\DF X$$
on the full flag manifold 
(which in view of the inclusion $\DF X\subset\ol X^{Fins}$ is a restriction of the action on the Finsler compactification)
and the restricted domains 
$$\OmF=\Om\cap\DF X= \DF X - \ThF (\Lat)$$
the following holds:
The action $\Ga\acts\OmF$ is {\em properly discontinuous} as long as $\Th$ is {\em fat} 
in the sense that $W=\Th\cup w_0\Th$,
and {\em cocompact} as long as $\Th$ is {\em slim} 
in the sense that $\Th\cap w_0\Th=\emptyset$,
compare \cite[Thms.\ 1.5 and 1.8]{coco15}. 
A thickening is {\em balanced} iff it is both fat and slim. 

These results are reminiscent of, and were motivated by 
the notion of {\em stability} and the construction of the Mumford quotient 
for actions $G\acts V$ of semisimple algebraic groups on projective varieties 
in {\em Geometric Invariant Theory}.
Our domains of proper discontinuity correspond in GIT to the stable part of the action,
and the balancedness
to the case when ``semistable$\Ra$stable". 

There is not only an analogy, but also a concrete relation
between our theory and Mumford's GIT. 
For instance, in the case of certain configuration spaces
we recover the GIT quotient, 
cf.\ \cite[Ex.\ 3.38 and \S 7.4]{coco15}.

Our proper discontinuity and cocompactness 
results extend the scope of the earlier work in \cite{GW}, 
where domains of proper discontinuity and cocompactness were constructed 
for actions 
of Anosov subgroups of classical semisimple Lie groups on flag manifolds, and
in the case of Anosov subgroups of general semisimple Lie groups 
only for actions 
on the quotients $G/AN$
which are fiber bundles over the full flag manifold $G/B$ with compact fibers.

{\bf 5.\ Related work.} 
The results discussed in this paper 
imply most of the main results of the papers 
\cite{GGKWmproper, GGKWtame}
which appeared later.
We also note that 
the argument for the most general result in \cite{GGKWtame}
is still incomplete as of now.

{\bf 6.\ Acknowledgements.} 
The first author was supported by the NSF grant DMS-12-05312. 
He also thanks the Korea Institute for Advanced Study  for its hospitality. 
The last author was supported by the grant Mineco MTM2012-34834.

\noindent M.K.: Department of Mathematics, 
University of California, Davis, 
CA 95616, USA\\
email: kapovich@math.ucdavis.edu

\noindent B.L.: Mathematisches Institut,
Universit\"at M\"unchen, 
Theresienstr. 39, 
D-80333 M\"unchen, Germany, 
email: b.l@lmu.de

\noindent J.P.: Departament de Matem\`atiques, 
Universitat Aut\`onoma de Barcelona, 
E-08193 Bellaterra, Spain, 
email: porti@mat.uab.cat

\end{document}